# Fredrickson–Andersen Model and Noisy Majority Vote Process on Nonamenable Graphs

Damiano De Gaspari*

August 9, 2025

**Abstract**

We study the Fredrickson-Andersen $j$-spin facilitated model and the noisy majority vote process on connected infinite graphs satisfying suitable expansion properties. For the former, we consider the out-of-equilibrium regime where the density of facilitating sites is close to 1, both for the equilibrium product measure and for the initial configuration, and we show exponential convergence to equilibrium. For the latter, we prove the existence of multiple equilibrium measures, generalising recent results by J. Ding and F. Huang (2025). Our proofs build on the framework of decorated set systems introduced by I. Hartarsky and F. Toninelli (2024) and establish exponentially decaying tails for the diameter of the space-time cluster of zeros containing a fixed vertex for both perturbed bootstrap percolation and consensus processes. The results are essentially sharp on hyperbolic lattices and we further show how some of the borderline cases can be handled using Toom contours, in the reformulation by J. M. Swart, R. Szabó and C. Toninelli (2022).

## 1 Introduction

Kinetically constrained models (KCMs) are a class of interacting particle systems introduced in the statistical physics literature [FA84; FA85] to investigate the *liquid-glass transition.* Typically defined on discrete lattices, these models evolve according to a stochastic Glauber dynamics in which each spin is resampled from the equilibrium Bernoulli($q$) product measure only if a local constraint is satisfied. The parameter $q \in [0,1]$ represents the density of facilitating sites and, crucially, the constraint does not involve the spin being updated. KCMs are thus designed to capture glassy behaviour without relying on non-trivial equilibrium properties, reflecting the viewpoint that glassy physics is a fundamentally dynamical phenomenon [RS03; BG13]. The underlying physical intuition that the motion of a molecule in a dense liquid can be hindered by the presence of too many nearby molecules is encoded mathematically by the fact that the jump rates can be degenerate. This results in a non-attractive, cooperative and heterogeneous dynamics, which exhibits multiple invariant measures and possibly ergodicity-breaking phase transitions.

Most of the work on KCMs has been conducted on the $d$-dimensional Euclidean lattice, considering general local constraints and focusing on the equilibrium regime, namely when the initial configuration is sampled from the equilibrium Bernoulli product measure. In this setting, many results of different kinds have been obtained and we refer to [HT24a], [Har22, Chapter 1] and references therein for a comprehensive overview. In contrast, the out-of-equilibrium behaviour of KCMs is much less understood. Physically, it corresponds to the analysis of the slow dynamics following a sudden quench from one temperature to another, see [MS07; LMSBG07]. Rigorous investigation began with [CMST10] and the most general results to date are those of [HT24b], whose Section 4 surveys more precise statements for specific models, see also [HT24a, Chapter 7].

---

*Technische Universität Wien, Institut für Stochastik und Wirtschaftsmathematik, Wiedner Hauptstraße 8-10, A-1040, Vienna, Austria. E-mail: `damiano.de.gaspari@tuwien.ac.at`

On general graphs, the study of arbitrary local constraints becomes less central and attention is typically focused on the $j$-spin facilitated Fredrickson-Andersen model, where a spin may be updated only if at least $j$ of its nearest neighbours are in the facilitating state. The works [CLR79; SBT05; SLC06; PR24] highlight interest from the physics community. On the mathematical side, when the underlying graph is a regular tree, [MT13; CMRT15] provide a comprehensive picture by identifying a threshold for ergodicity, proving exponential convergence to equilibrium in the ergodic regime and power-law scaling at criticality. In [CMRT09] general infinite connected graphs are considered, but only for $j = 1$, which significantly simplifies the analysis, as in this case the dynamics is non-cooperative. In the present work, we establish relationships between the parameter $j$ and various notions of expansion of the underlying graph that guarantee exponential convergence to equilibrium, in the sense of local functions, in the regime where the density of facilitating sites is close to 1, both for the equilibrium product measure and for the initial configuration. The result mirrors what was obtained in [HT24b] for the Euclidean lattice. In the specific case of hyperbolic lattices, our conditions encompass essentially all instances in which ergodicity is expected according to [STBT10].

The majority vote process with noise $\varepsilon \in [0, 1]$ is an interacting particle system describing the evolution of two opinions in a population whose social network is represented by the underlying graph. Each time a spin is updated, with probability $1 - \varepsilon$ it is set equal to the majority of the opinions in its voting neighbourhood, while with probability $\varepsilon$ it is resampled uniformly at random. If $\varepsilon$ is sufficiently large, [BG21, Proposition 1] shows that the dynamics converges exponentially fast to a unique equilibrium for every initial configuration (see also [Lig05, Theorem I.4.1]). In contrast, much less is understood about the behavior of the system in the low-noise regime. On the one-dimensional Euclidean lattice there exists a unique equilibrium for all $\varepsilon > 0$, see [Gra82], [Lig05, Example III.2.12] and [LMS90, Examples 1.1, 6.1]. However, the equilibrium behaviour on higher-dimesional Euclidean lattices remains a long-standing open problem and essentially nothing is known in this setting (see [Lig05, Sections I.9, III.8] for some open problems). On the infinite $d$-regular tree, [BG21] showed that there exist uncountably many mutually singular equilibria if $d \geq 5$ and the recent [DH25], employing different techniques, extended the non-uniqueness result to $d \geq 3$. Finally, [BCOTT16] investigates the regime $\varepsilon = 0$ on unimodular transitive graphs, showing that each spin either stabilises or oscillates with period two. Within the physics literature, [VPTMA20] examines the cases of mean-field behaviour and random underlying graph. In the present work, we establish the existence of multiple invariant measures when the graph expansion constant is sufficiently large and $\varepsilon$ is sufficiently small. The setting of hyperbolic lattices again serves as a natural concrete example and has also attracted recent interest in the physics literature, see for instance [WH10; ECDF19]. Moreover, for $d$-regular trees, we show how Toom contours can be used to sharpen our analysis, providing an alternative derivation of the results in [DH25].

The decision to present these distinct models and results within a single work is motivated by the common tools and proof strategies we employ. In Section 2, after introducing the necessary notation, we give precise definitions of the two models described, as well as of bootstrap percolation and consensus processes. These will serve as tools in our analysis, but are also of independent interest. We then state our main results and provide further comments. Section 3 is devoted to applying our results to the setting of hyperbolic lattices. Our main original contribution is Section 4, where we construct graphs that explain how a certain spin value propagates through space-time to reach an arbitrary given point. In particular, we show that these constructions satisfy the combinatorial properties required to fit within the framework of decorated set systems of [HT24b], which we recall in Section 5 and which serve as the workhorse of our approach. Finally, Section 6 collects the proofs of the results stated in Section 2.

## 2 Models and Results

We work on $G = (V, E)$, an unoriented connected infinite graph without self-loops. For each $x, y \in V$, we write $x \sim y$ if $\{x, y\} \in E$ and we denote by $N(x) \coloneqq \{y \in V \mid y \sim x\}$ the set of the neighbours of $x$ in $G$. We further set $\deg(x) \coloneqq |N(x)|$, where $|\cdot|$ denotes cardinality and we assume that $\Delta(G) \coloneqq \sup\{\deg(x) \mid x \in V\} < \infty$, which is sufficient to ensure that the interacting particle systems introduced below are well-defined, see e.g. [Lig99; Lig05]. Given $G$, we denote by $\Omega \coloneqq \{0, 1\}^V$ the set of configurations $(\eta_x)_{x \in V}$ that assign either 0 or 1 to each vertex of the graph and we use the symbols **0** and **1** for those that are everywhere zero and one, respectively. Moreover, if $\eta$ is a Markov process with state space $\Omega$ and càdlàg paths, as will always be the case below, we denote by $\eta_x(t)$ the value of the process at vertex $x \in V$ at time $t \in [0, \infty)$ and by $\eta_x(t^-)$ the limit $\lim_{s \to t^-} \eta_x(s)$. Let $\mu$ be a generic probability measure on $\Omega$. We write $\eta^\mu$ and $\eta^\xi$, for $\xi \in \Omega$, to say that the configuration at time 0 has been sampled according to the measures $\mu$ and $\delta_\xi$, respectively. For any given $x \in V$, we use the function $c_x \colon \Omega \to \mathbb{N}$, defined by $\eta \mapsto c_x(\eta) \coloneqq \sum_{y \sim x} \eta_y$, to count how many neighbours of $x$ have value 1 in every given configuration (our convention is that $\mathbb{N}$ denotes the set of nonnegative integers). Finally, we recall that a function $f \colon \Omega \to \mathbb{R}$ is called *local* if there exist a finite set $\text{supp}(f) \subset V$ and a function $\tilde{f} \colon \{0, 1\}^{\text{supp}(f)} \to \mathbb{R}$ such that $f(\eta) = \tilde{f}(\eta|_{\text{supp}(f)})$ for every $\eta \in \Omega$.

To define the dynamics of the models under consideration, we introduce the collections of independent random variables $(P_x)_{x \in V}$, $(U_x(t))_{x \in V, t \in P_x}$ and $(L_x(t))_{x \in V, t \in \mathbb{N}_{\geq 1}}$ with the following respective distributions: Poisson point process on $[0, \infty)$ of intensity 1, uniform on $[0, 1]$ and again uniform on $[0, 1]$. For convenience, we assume that these random variables are defined on the same probability space, with probability measure and expectation denoted by $\mathbb{P}$ and $\mathbb{E}$. As usual, when considering random initial configurations, we implicitly extend the probability space so that the initial configurations are independent of the dynamics. With a slight abuse of notation, we continue to use the symbols $\mathbb{P}$ and $\mathbb{E}$ even when averaging also over the initial configuration. Finally, for each fixed $\varepsilon \geq 0$, we say that $(x, t) \in V \times \mathbb{N}_{\geq 1}$ is an *update point* if $L_x(t) \geq \varepsilon$ and a *noise point* otherwise. This is consistent with [BG21], whereas in [SST25b] and [HT24b] noise points are referred to as defective sites and deaths, respectively.

Throughout the remainder of this work, unless stated otherwise, $j$ always denotes an integer satisfying $1 \leq j \leq \delta(G) \coloneqq \inf\{\deg x \mid x \in V\}$. We are now ready to define the two processes we are mainly interested in.

### 2.1 Fredrickson-Andersen model and Bootstrap Percolation

**Definition 2.1** (FA-$j$f)**.** We call *Fredrickson-Andersen j-spin facilitated model on $G$ with equilibrium parameter $q \in [0, 1]$* a continous-time Markov process with state space $\Omega$, càdlàg paths, dynamics defined by

$$\eta_x(t) = \begin{cases} \mathbb{1}_{\{U_x(t) \leq q\}} & \text{if } t \in P_x \text{ and } c_x(\eta(t^-)) \geq j, \\ \eta_x(t^-) & \text{otherwise,} \end{cases} \tag{2.1}$$

for every $x \in V$ and $t \in \mathbb{R}_{>0}$, and initial configuration $\omega(0)$ chosen according to a prescribed probability distribution $\mu$ on $\Omega$.

**Definition 2.2** (BP with $\varepsilon$-noise)**.** We call *$j$-neighbour Bootstrap Percolation with $\varepsilon$-noise on $G$* a discrete-time Markov process with state space $\Omega$, dynamics inductively defined by

$$\omega_x(t) = \begin{cases} \max\{\omega_x(t-1), \mathbb{1}\{c_x(\omega(t-1)) \geq j\}\} & \text{if } L_x(t) \in [\varepsilon, 1], \\ 0 & \text{if } L_x(t) \in [0, \varepsilon), \end{cases} \quad (2.2)$$

for every $x \in V$ and $t \in \mathbb{N}_{\geq 1}$, and initial configuration $\omega(0)$ chosen according to a prescribed probability distribution $\mu$ on $\Omega$. If $\varepsilon = 0$, we simply call it Bootstrap Percolation.

For any given initial condition, almost surely there exists a unique process $\eta$ satisfying Definition 2.1. This follows from the general graphical construction framework for interacting particle systems, originally introduced in [Har78]. Common references include [Swa22] and [Lig05, Section III.6]. Analogous existence and uniqueness results hold for the other processes defined below. The generator $\mathcal{L}$ of $\eta$ can also be constructed in a standard way in terms of its action on local functions. It is a non-positive self-adjoint operator on $L^2(\Omega, \mu_q)$ with domain $\text{Dom}(\mathcal{L})$, [Lig05, Sections I.3, IV.4]. We remark that $\eta$ is not attractive, in the terminology of [Lig05, Section III.2] (or monotone, in the terminology of [Swa22]), i.e. the dynamics (2.1) does not preserve the stochastic (partial) order. It is instead straightfoward to verify that the dynamics with equilibrium parameter $q$ satisfies the detailed balanced condition with respect to the Bernoulli product measure $\mu_q := \text{Ber}(q)^{\otimes V}$, as for any configuration $\xi \in \Omega$ the value $c_x(\xi)$ is independent of $\xi_x$. Therefore, $\eta$ is reversible with respect to $\mu_q$, [Lig05, Section IV.2]. However, due to the fact that the jump rates of the process are not bounded away from zero, irreducibility is not guaranteed. Indeed, $\mu_q$ is not the only invariant measure: $\delta_{\mathbf{0}}$ is also clearly invariant, as $c(\mathbf{0}) = \mathbf{0}$ implies that no vertex is ever updated. It is thus natural to define the *ergodicity critical parameter*

$$q_c := \inf\{q \in [0, 1] \mid 0 \text{ is a simple eigenvalue of } \mathcal{L}\}.$$

It turns out that Boostrap Percolation $\omega$ (without noise) is closely related to the Fredrickson-Andersen model and to KCMs in general, see for example [HT24a, Chapter 3]. In particular, it holds that

$$q_c = q_c^{\text{BP}} := \inf\left\{q \in [0, 1] \,\bigg|\, \lim_{t \to \infty} \omega^{\mu_q}(t) = \mathbf{1}\right\} \quad (2.3)$$

and moreover, for any $q > q_c$, the eigenvalue 0 is simple for $\mathcal{L}$ (the pointwise limit of $\omega(t)$ as $t \to \infty$ exists because the trajectories of $\omega$ are monotone). On Euclidean lattices, this is the content of [CMRT08, Proposition 2.5]. However, as already noted in [CMRT09, Remark 4], the proof applies to any graph of bounded degree, and so also to our setting. Determining whether $q_c^{\text{BP}}$ is nontrivial, namely whether $q_c^{\text{BP}} \in (0, 1)$, has been the subject of extensive study. On $d$-dimensional Euclidean lattices, for the specific update rule given in (2.2), one has $q_c^{\text{BP}} = \mathbb{1}\{j > d\}$, as shown in [Sch92]. A complete classification of arbitrary local update rules has only recently been achieved, see [BSU15; BBPS16; BBMS22; BBMS24]. On nonamenable graphs, instead, there exist values of $j$ for which (2.2) defines models with a nontrivial critical density $q_c^{\text{BP}} \in (0, 1)$, as shown in [BPP06; FS08; STBT10]. This is the setting of the present work.

In its version with noise, Bootstrap Percolation has been studied as a perturbed monotone cellular automaton. Because of stochastic monotonicity, the sequence of probability measures $(\mathbb{P}[\omega^{\mathbf{1}}(t) \in \cdot])_{t \in \mathbb{N}}$ converges weakly to the *upper invariant measure* $\bar{\nu}$ of the process and it makes sense to define the *stability threshold*

$$\varepsilon_c^{\text{BP}} := \sup\{\varepsilon \in [0, 1] \mid \bar{\nu} \neq \delta_{\mathbf{0}}\}. \quad (2.4)$$

The work of [Too80], recently revisited and extended in [SST25a; SST25b], investigates necessary and sufficient conditions so that $\varepsilon_c^{\text{BP}} > 0$. On Euclidean lattices, this is shown to be equivalent to the property that any initial configuration with only finitely many vertices in state 0 eventually evolves to **1** under the dynamics without noise. On general graphs, the situation is less clear. As a byproduct of our arguments, we obtain that $\varepsilon_c^{\text{BP}} > 0$ whenever $j$ is $\omega$-good in the sense of Definition 2.5 below, see Corollary 2.10.

### 2.2 Noisy Majority Vote and Consensus Processes

**Definition 2.3** ($NMVP$). We call $\varepsilon$-*noise Majority Vote Process on* $G$ a discrete-time Markov process with state space $\Omega$ and dynamics inductively defined by

$$\sigma_x(t) = \begin{cases} 0, & \text{if } L_x(t) \in [\varepsilon, 1] \text{ and } c_x(\sigma(t-1)) < \frac{\deg(x)}{2}, \\ 1, & \text{if } L_x(t) \in [\varepsilon, 1] \text{ and } c_x(\sigma(t-1)) > \frac{\deg(x)}{2}, \\ \sigma(t-1) & \text{if } L_x(t) \in [\varepsilon, 1] \text{ and } c_x(\sigma(t-1)) = \frac{\deg(x)}{2}, \\ 0 & \text{if } L_x(t) \in [0, \frac{\varepsilon}{2}), \\ 1 & \text{if } L_x(t) \in [\frac{\varepsilon}{2}, \varepsilon), \end{cases}$$

for every $t \in \mathbb{N}_{\geq 1}$.

**Definition 2.4** ($j$-$CP$). We call $j$-*threshold Consensus Process with* $\varepsilon$-*noise on* $G$ a discrete-time Markov process with state space $\Omega$, dynamics inductively defined by

$$\chi_x(t) = \begin{cases} \mathbb{1}\{c_x(\chi(t-1)) \geq j\} & \text{if } L_x(t) \in [\varepsilon, 1], \\ 0 & \text{if } L_x(t) \in [0, \varepsilon), \end{cases} \tag{2.5}$$

for every $t \geq 1$, and initial configuration $\chi(0)$ chosen according to a prescribed probability distribution $\mu$ on $\Omega$.

The definition of the Noisy Majority Vote Process given above matches that in [DH25], while [BG21] considers its continuous-time analogue, in which updates occur at the time points of the Poisson point processes $(P_x)_{x \in V}$. Moreover, our consensus processes correspond to discrete-time analogues of the threshold processes of [BG21], but in our case they are biased towards state 0 instead of state 1. It is important to note the difference between the dynamics (2.4) and that of Bootstrap Percolation defined in (2.2): in consensus processes, unlike in Bootstrap Percolation, a vertex can change its value from 1 to 0 even in the absence of noise. Finally, since the dynamics of $\chi$ is stochastically monotone, it is natural to consider $\bar{\nu}$, the upper invariant measure of the process $\chi$, and define the *stability threshold*

$$\varepsilon_c^{\text{CP}} := \sup\{\varepsilon \in [0, 1] \mid \bar{\nu} \neq \delta_{\mathbf{0}}\}, \tag{2.6}$$

analogoulsy to what already done above for Boostrap Percolation.

### 2.3 Results

We start by describing the different notions of "expansion" that will be used in the results below. First, we define the *edge expansion constant of* $G$ as

$$\Phi_E(G) := \inf\left\{ \frac{|\partial_E K|}{|K|} \,\bigg|\, \emptyset \neq K \subset V, \ K \text{ finite} \right\}, \tag{2.7}$$

where $\partial_E K \subset E$ is the set of edges of $G$ which have exactly one endpoint in $K$. Similarly, the *vertex expansion constant of* $G$, denoted by $\Phi_V(G)$, is defined by replacing $\partial_E K$ in (2.3) with $\partial_V K := \{x \notin K \mid \exists y \in K \ x \sim y\}$ (see for example [LP17, Chapter 6]). Clearly, one always has $\Phi_V \leq \Phi_E$. Graphs for which either of these constants is strictly positive are called *nonamenable*.

For every natural number $k$ and vertex $x \in V$, we denote by $B_k(x) := \{y \in G | \text{dist}(x,y) \leq k\}$ and $S_k(x) := \{y \in G \mid \text{dist}(x,y) = k\}$ the ball and the sphere, respectively, of radius $k$ centered at $x$ with respect to the graph distance. We define

$$\bar{j} := \max\{j \in \mathbb{N} \mid \exists\, r \in V, \forall k \in \mathbb{N}, \forall x \in S_k(r), j \leq |S_{k+1}(r) \cap N(x)|\} \tag{2.8}$$

We observe that the maximum is taken over a nonempty set, since 0 trivially belongs to the set on the right-hand side. Moreover, the hypothesis that the degrees of $G$ are uniformly bounded guarantees that $\bar{j}$ is finite.

Those notions allow us to formulate several, generally incomparable, assumptions on the parameter $j$ under which our results hold. For brevity, we denote by $\Xi$ either the $j$-threshold Consensus Process $\chi$ or the $j$-neighbour Bootstrap Percolation $\omega$ with $\varepsilon$-noise on the graph $G$. We adopt this convention for the remainder of the work. Moreover, we give the following definition in order to summarize the various conditions we assume on $j$ in our results.

**Definition 2.5.** Let $j$ be an integer such that $1 \leq j \leq \delta(G)$.
We say that $j$ is $\omega$-*good* if it satisfies at least one of the following conditions:
(a) $j < \Phi_E + 1$; (b) $j \leq \bar{j}$, where $\bar{j}$ is given by (2.3).
We say that $j$ is $\chi$-*good* if it satisfies at least one of the following conditions:
(i) $2j < 3\Phi_E - \Delta + 2$; (ii) $G$ is bipartite and $j < \Phi_E + 1$; (iii) $j < \Phi_V$; (iv) $G$ is the $d$-regular tree, with $d \geq 2$, and $j < d$.

We are now ready to state our two main results.

**Theorem 2.6** (FA-$j$f – exponential convergence to equilibrium)**.** *Let $(\eta^{\mu_p}(t))_{t \geq 0}$ be the Fredrickson-Andersen $j$-spin facilitated model on $G$ with initial distribution $\mu_p$, equilibrium parameter $q$ and $\omega$-good $j$. Then there exist $\varepsilon, c > 0$ such that for every $p, q \in [1 - \varepsilon, 1]$, for every local function $f \colon \Omega \to \mathbb{R}$ and for every $t \in [0, \infty)$*

$$|\mathbb{E}[f(\eta^{\mu_p}(t))] - \mu_q(f)| \leq c^{-1} |\mathrm{supp} f|\, \|f\|_\infty\, e^{-ct}, \tag{2.9}$$

*where $\mu_q(f)$ denotes the integral of $f$ with respect to the measure $\mu_q$.*

In fact, as already noted in [HT24b], it is clear from the proof that the result remains valid for any initial distribution $\mu$ that stochastically dominates $\mu_p$. In light of the discussion surrounding (2.1), we conjecture that Theorem 2.6 holds for all $p, q \in (q_c, 1]$. Indeed, the analogous result in the setting of Euclidean lattices established in [HT24b] allows for $p$ in the whole regime where the spectral gap of the generator of $\eta$ is positive. This is obtained thanks to a renormalization in space, an approach which does not immediately apply to general graphs. Extending convergence to equilibrium beyond the perturbative regime for the equilibrium measure appears even more challenging and has been achieved only on very specific models, see [HT24a, Chapter 7]. In the setting of this paper, the only case we expect to be within reach using the current tools is $j = 1$.

**Theorem 2.7** (NMVP and $j$-CP – non-ergodicity)**.** *Let $(\chi(t))_{t\geq 0}$ and $(\sigma(t))_{t\geq 0}$ be the $j$-threshold Consensus Process and the Majority Vote Process with $\varepsilon$-noise on $G$, respectively. Assume that $j$, respectively $\lfloor \Delta/2 + 1 \rfloor$, is $\chi$-good and that $\varepsilon$ is sufficiently small. Then $\chi$, respectively $\sigma$, has multiple invariant measures and there exists $0 < c < \tfrac{1}{2}$ such that for every $t \in \mathbb{N}_{\geq 1}$ and for every $x \in V$*

$$\mathbb{P}[\chi_x^{\mathbf{1}}(t) = 0] \leq c, \qquad \mathbb{P}[\sigma_x^{\mathbf{1}}(t) = 0] \leq c.$$

*Remark* 2.8 (Regular trees)**.** We assume throughout this remark that the underlying graph is the $d$-regular tree. Consider the Noisy Majority Vote Process. Item (ii) of Definition 2.5 holds for $d \geq 5$ and in this case we establish Theorem 2.7 using the history graphs introduced in Section 4.1, thereby providing a new proof of one of the main results in [BG21]. On the other hand, item (iv) of Definition 2.5 holds for all $d \geq 3$ and in this case we derive the result of Theorem 2.7 via the Toom contour argument described in Section 4.2. This yields an alternative proof of the results in [DH25] and confirms the prediction, stated in the introduction of [BG21], that the result should be obtained using some version of Toom's construction. Finally, consider the $j$-threshold Consensus Process. If $j = d$, the result is clearly false, as two neighbours whose values are simultaneously set to zero by the noise remain zero forever. In this sense, item (iv) is sharp. We observe that it also applies to the 1-threshold Consensus Process on the 2-regular tree, in which case Theorem 2.7 reduces to the classical fact that the critical probability for the existence of an infinite connected component in oriented site percolation on $\mathbb{Z}^2$ is bounded away from 1.

In the proofs of Theorems 2.6 and 2.7 we also recover the following proposition, which is of independent interest. We recall that the *strong product* of two graphs $G_i = (V_i, E_i)$, $i = 1, 2$, is the graph $G_1 \boxtimes G_2$ whose vertices are the elements of $V_1 \times V_2$ and whose edges are of the form $\{(x_1, x_2), (y_1, y_2)\}$ and satisfy one of the following conditions:

- $x_1 = y_1$ and $\{x_2, y_2\} \in E_2$;
- $\{x_1, y_1\} \in E_1$ and $x_2 = y_2$;
- $\{x_1, y_1\} \in E_1$ and $\{x_2, y_2\} \in E_2$.

The graph defined in the same way, but allowing only the edges given by the first two items above is called the *Cartesian product* of $G_1$ and $G_2$, denoted by $G_1 \square G_2$. Although $G_1 \square G_2$ is the natural choice to model the space-time structure of our processes, the stronger connectivity of $G_1 \boxtimes G_2$ is more appropriate for the construction of history graphs introduced in the next section and yields stronger results.

**Proposition 2.9.** *Let $\Xi$ be either the $j$-threshold Consensus Process $\chi$ or the $j$-neighbour Bootstrap Percolation $\omega$ with $\varepsilon$-noise on $G$. Let $(x, t)$ be a fixed vertex of $G \boxtimes \mathbb{N}$ and denote by $\mathscr{C}_{x,t}$ the connected component in $G \boxtimes \mathbb{N}$ of the set $\{(y, s) \in V \times \mathbb{N} \mid \Xi^1_y(s) = 0\}$ containing $(x, t)$. Assume that $j$ is $\Xi$-good. Then there exist $c > 0$ and $\widetilde{\varepsilon} > 0$ such that for every $\varepsilon \in [0, \widetilde{\varepsilon}]$ and integer $\ell \geq 0$*

$$\mathbb{P}[\operatorname{diam}(\mathscr{C}_{x,t}) \geq \ell] \leq \varepsilon^{c(\ell+1)}, \tag{2.10}$$

*with the convention that $\operatorname{diam}(\emptyset) = -\infty$.*

**Corollary 2.10.** *Let $\varepsilon_c^{\mathrm{BP}}$ and $\varepsilon_c^{\mathrm{CP}}$ denote the stability thresholds of the $j$-neighbour Bootstrap Percolation and the $j$-threshold Contact Process with $\varepsilon$-noise on $G$, respectively, as defined in (2.1) and (2.2). If $j$ is $\omega$-good, then $\varepsilon_c^{\mathrm{BP}} > 0$. If $j$ is $\chi$-good, then $\varepsilon_c^{\mathrm{CP}} > 0$.*

## 3 Hyperbolic lattices: a case study

In this section, we consider hyperbolic lattices, which offer a concrete class of graphs, other than regular trees, on which to test the strength of the assumptions of Theorems 2.6 and 2.7, respectively. In other words, we aim to determine whether $j$ is $\omega$-good or $\chi$-good, respectively, in the sense of Definition 2.5, in all cases where the results are expected to hold.

Let $d$ and $f$ be integers such that $d, f \geq 3$ and $(d-2)(f-2) > 4$. We define the *hyperbolic lattice* $H(d, f)$ to be the graph obtained as the tessellation of the Poincaré disc consisting of regular $f$-gons with interior angles equal to $2\pi/d$. In other words, $H(d, f)$ is the regular plane graph of degree $d$ whose dual is regular and of degree $f$. A derivation of these tessellations in terms of the so called Fuchsian groups, discrete subgroups of the group of isometries of the hyperbolic plane, can be found, for example, in [CM80; Kat92; Ive92]. We remark that, although either $H(d,f)$ or $H(f,d)$ is a Cayley graph, as shown in [CK96], not all hyperbolic lattices are. In particular, those that are not fall outside of the scope of [SST25a].

The edge expansion constant (2.3) for hyperbolic lattices has been computed explicitly, with two independent proofs given in [HJL02; HS03]:

$$\Phi_E(H(d,f)) = (d-2)\sqrt{1 - \frac{4}{(d-2)(f-2)}}.$$

To the best of the author's knowledge, no explicit formula is known in general for the vertex expansion constant. However, good lower bounds are provided by [HP21, Theorems 4, 5], which establish that for every plane graph $G_1$ with minimum degree $d \geq 7$ and $G_2$ with minimum degree $d \geq 5$ and triangle-free, the following inequalities hold:

$$\Phi_V(G_1) \geq \alpha_d, \qquad \Phi_V(G_2) \geq \alpha_{d+2}, \qquad \alpha_d := \frac{d - 6 + \sqrt{(d-2)(d-6)}}{2}, \tag{3.1}$$

with the hyperbolic lattices $G_1 = H(d, 3)$ and $G_2 = H(d, 4)$ realizing equality.

**Bootstrap Percolation with noise** Let us consider $j$-neighbour Bootstrap Percolation (i.e. $\varepsilon = 0$) on $H(d,f)$ with initial distribution $\mu_q$. As shown in [STBT10, Appendix], if $j$ does not satisfy

$$1 \leq j \leq d-3 \text{ if } f = 3, \qquad 1 \leq j \leq d-2 \text{ if } f \geq 4, \tag{3.2}$$

then, with positive probability, there exist finite clusters of zeros in the initial configuration which remain zero for all future times. Therefore, (2.9) can only be true if $j$ is in the above ranges. It is not difficult to verify that item (a) of Definition 2.5 is true in all but the following cases:

$$f = 4, d \geq 5, j = d-2, \qquad f = 3, d \geq 7, j = d-3.$$

The case $f = 4$ is covered by item (b) of Definition 2.5. Indeed, let $x \in V$ and $k \in \mathbb{N}_{\geq 1}$ be arbitrary. By the Construction Algorithm of [BD10, Section 2], each vertex $y \in S_k(x)$ is of one of two types: either adjacent to exactly one vertex in $B_k(x)$, or adjacent to exactly two vertices in $B_k(x)$ (at most one of which is in $S_k(x)$). In both cases, $|S_{k+1}(x) \cap N(y)| \geq d-2$.

The case $f = 3$, $d \geq 7$, $j = d-3$, instead, does not satisfy our definition of $\omega$-good and thus remains open.

**Consensus and Noisy Majority Vote Process** We first consider the Noisy Majority Vote Process. We seek to determine for which values of $d$ the number $\lfloor d/2 + 1 \rfloor$ is $\chi$-good when the underlying graph is the hyperbolic lattice $H(d,f)$. Using the first inequality in (3), it is immediate to verify item (iii) of Definition 2.5 for every $H(d,f)$ with $d \geq 12$. Using the second inequality in (3), we find that this is also the case for every $H(d,f)$ with $d \geq 8$ and $f \geq 4$. The cases $d = 7, 8$ satisfy item (i) if $f \geq 6$. Finally, recall that every planar graph in which every face is bounded by an even number of edges is bipartite. We deduce that item (ii) is satisfied by any $H(d,f)$ such that $f$ is even and one of the following is satisfied: $d \geq 7$ and $f \geq 4$; $d = 6$ and $f \geq 6$; $d = 5$ and $f \geq 8$. As already mentioned in the introduction, to the best of our knowledge, there is no nontrivial graph, apart from the one-dimensional Euclidean lattice, for which it is known (or even conjectured) that the process has a unique invariant measure for sufficiently small $\varepsilon$.

We now consider the $j$-threshold Consensus Process with $\varepsilon$-noise $\chi$ on $H(d,f)$. First of all, we observe that Theorem 2.7 cannot hold for all admissible values of $d$, $f$ and $j$. Indeed, suppose that at some time $t$ the configuration $\chi(t)$ is zero on some finite subset of the vertices. Then this subset remains zero for all future times if and only if the same is true for $j$-neighbour Bootstrap Percolation. Therefore, if $j$ does not satisfy (3), then for every $\varepsilon > 0$ and every initial configuration $\xi$ we have that $\lim_{t\to\infty} \chi^\xi(t) = \mathbf{0}$ almost surely. In particular, $\delta_{\mathbf{0}}$ is the unique invariant measure of $\chi$. On the other hand, an analysis analogous to that carried out above for the Noisy Majority Vote Process shows that $j$ is $\chi$-good for nearly all the cases described by (3). We conjecture that Theorem 2.7 also holds in the remaining ones.

## 4 Leveraging expansion

In this section we consider history graphs and Toom cycles and contours, deterministic objects which provide information about the event that the process assigns value 0 to a given vertex at a given time. We also establish several of their combinatorial properties, which are crucial for the effective application of the decorated set system framework of Section 5. Compared to history graphs, the construction of Toom cycles and contours is more sophisticated and in the absence of structure in the underlying graph, it is generally unclear how to derive the aforementioned combinatorial properties. However, it can yield sharper results.

### 4.1 History graphs

We work on the oriented space-time graph $\mathscr{G} := (\mathscr{V}, \mathscr{E})$, defined by $\mathscr{V} := V \times \mathbb{N}$ and $\mathscr{E} := \mathscr{S} \sqcup \mathscr{O}$, where $\sqcup$ denotes disjoint union, $\mathscr{S} := \{((x,t),(x,t-1)) \in (V \times \mathbb{N})^2 \mid x \in V\}$ is the set of *straight* edges and $\mathscr{O} := \{((x,t),(y,t-1)) \in (V \times \mathbb{N}_0)^2 \mid x,y \in V, x \sim y\}$ is the one of *oblique* edges. In words, $\mathscr{E}$ is the set of all possible oriented edges connecting a vertex $x$ at time $t$ to either itself (straight edge) or one of its neighbours (oblique edge) at time $t-1$. Finally, given $t \in \mathbb{N}$ and $\mathscr{U} \subset \mathscr{V}$, we denote by $U_t := \{x \in V \mid (x,t) \in \mathscr{U}\}$ the space projection of the $t$-section of $\mathscr{U}$. The general convention is that we use the main font for sets of vertices belonging to the base graph and the math script font for sets of vertices belonging to the space-time graph.

**Definition 4.1** (History graph). Let $j$ be an integer satisfying $1 \leq j \leq \delta(G)$ and let $(o,T) \in \mathscr{V}$ be an arbitrary space-time vertex. Consider a finite, oriented, connected subgraph $(\mathscr{U}, \mathscr{F})$ of $\mathscr{G}$ and a subset $\mathscr{U}^* \subset \mathscr{U}$. We say that $\mathscr{H} = (\mathscr{U}, \mathscr{F}, \mathscr{U}^*)$ is a *history graph* for $(o,T)$ for the $j$-threshold Consensus Process on $G$ if the following properties are satisfied:

(1) $(o,T) \in \mathscr{U}$ and all other vertices in $\mathscr{U}$ have time coordinate in $[1,T]$;

(2) all the vertices of $\mathscr{U}$ with time coordinate 1 are in $\mathscr{U}^*$;

(3) no vertex in $\mathscr{U}^*$ has outgoing edges;

(4) each vertex $x \in \mathscr{U} \setminus \mathscr{U}^*$ has exactly $\deg_G(x) - j + 1$ outgoing oblique edges.

We say that $\mathscr{H}$ is a *history graph* for $(o,T)$ for the $j$-neighbour Bootstrap Percolation on $G$ if it satisfies the above properties with (4) replaced by

(5) each vertex $x \in \mathscr{U} \setminus \mathscr{U}^*$ has exactly one outgoing straight edge and exactly $\deg_G(x) - j + 1$ outgoing oblique edges.

*Remark* 4.2. We observe that every history graph is uniquely determined by its edge set $\mathscr{F}$, since the subset $\mathscr{U}^*$ is then specified by item (3) of the definition. For this reason, it is not an abuse of notation to denote history graphs simply by $\mathscr{H} = (\mathscr{U}, \mathscr{F})$.

**Lemma 4.3.** *Let $\Xi^{\mathbf{1}}$ be either $j$-threshold Consensus Process or $j$-neighbour Bootstrap Percolation with $\varepsilon$-noise on $G$, started from the constant one configuration. Assume that $\Xi^{\mathbf{1}}_o(T) = 0$ for some $(o,T) \in \mathscr{V}$. Then there exists a history graph $\mathscr{H} = (\mathscr{U}, \mathscr{F}, \mathscr{U}^*)$ for $(o,T)$ for $\Xi$ which is present in $\Xi^{\mathbf{1}}$, which means that*

*(1) $\Xi^{\mathbf{1}}_x(t) = 0$ for every $(x,t) \in \mathscr{U}$;*

*(2) $(x,t)$ is a noise point, i.e. $L_x(t) < \varepsilon$, if and only if $(x,t) \in \mathscr{U}^*$.*

*Proof.* We construct $\mathscr{H}$ algorithmically. If $L_o(T) < \varepsilon$, we return $\mathscr{U} = \mathscr{U}^* = \{(o,T)\}$ and $\mathscr{F} = \emptyset$ and the algorithm ends. Clearly, this defines a history graph which satisfies the properties required by Lemma 4.3. Otherwise, if $L_o(T) \geq \varepsilon$, we initialize $\mathscr{U}_T = \bar{\mathscr{U}}_T = \{(o,T)\}$ and set $\mathscr{U}_T^* = \mathscr{F}_T = \emptyset$. We then iterate over $t$ in decreasing order, starting from $T$ and continuing until either $\mathscr{U}_t = \mathscr{U}_t^*$ or $t = 1$, performing the following three operations at each step:

(i) for each $(x,t) \in \bar{\mathscr{U}}_t$, identify a subset of $\{(y, t-1) \in \mathscr{V} \mid y \sim x, \Xi_y^{\mathbf{1}}(t-1) = 0\}$ of cardinality $\deg_G(x) - j + 1$;

(ii) if $\Xi$ is the Consensus Process, denote the subset identified in the previous step by $\mathscr{M}_t(x)$; if $\Xi$ is Bootstrap Percolation, define instead $\mathscr{M}_t(x)$ as union of that subset with $\{(x, t-1)\}$;

(iii) define

$$\mathscr{F}_t = \bigcup_{x \in \bar{\mathscr{U}}_t} \{((x,t),(y,t-1)) \mid (y,t-1) \in \mathscr{M}_t(x)\}, \quad \mathscr{U}_{t-1} = \bigcup_{x \in U_t} \mathscr{M}_t(x), \quad (4.1)$$
$$\mathscr{U}_{t-1}^* = \{y \in \mathscr{U}_{t-1} \mid L_y(t-1) < \varepsilon\}, \quad \bar{\mathscr{U}}_{t-1} = \mathscr{U}_{t-1} \setminus \mathscr{U}_{t-1}^*.$$

If the above iteration ends at a time $t > 1$, we further set $\mathscr{F}_s = \mathscr{U}_{s-1} = \emptyset$ for every $2 \leq s < t$. We return $\mathscr{U} = \bigcup_{t=1}^T \mathscr{U}_t$, $\mathscr{U}^* = \bigcup_{t=1}^T \mathscr{U}_t^*$ and $\mathscr{F} = \bigcup_{t=2}^T \mathscr{F}_t$ and the algorithm ends.

The only nontrivial step in the construction is the definition of $\mathscr{M}_t(x)$. Assume, for contradiction, that there exists $1 < t \leq T$ and a vertex $(x,t) \in \bar{\mathscr{U}}_t$ such that $|\{(y,t-1) \in \mathscr{V} \mid y \sim x, \Xi_y^{\mathbf{1}}(t-1) = 1\}| \geq j$. By the definitions in ((iii)) established at time $t+1$, we have $L_x(t) \geq \varepsilon$ and $\Xi_x^{\mathbf{1}}(t) = 0$. However, from $L_x(t) \geq \varepsilon$ it follows that the update rule of $\Xi$ applies (either (2.4) or (2.2)), which gives $\Xi_x^{\mathbf{1}}(t) = 1$, a contradiction. Thus, the algorithm is well-defined.

It follows immediately from the definitions that the output $\mathscr{H} = (\mathscr{U}, \mathscr{F}, \mathscr{U}^*)$ is a history graph satisfying the two properties stated in Lemma 4.3. In particular, item (2) of Definition 4.1 follows from the fact that $\Xi$ is started from the constant one configuration, so that the only way $\Xi_x(1) = 0$ is if $L_x(1) < \varepsilon$. □

*Remark* 4.4. In the terminology of [SST25b, Section 2.1], Bootstrap Percolation $\omega$ is a monotone cellular automaton, noise points are defective sites and $\mathscr{G}$ is the dependence graph associated to the collection of local maps which define $\omega$. Furthermore, (present) history graphs in $\omega^{\mathbf{1}}$ are very closely related to the so-called "(typed) explanation graphs", see [SST25b, Propositions 32, 33]. Elements of $\mathscr{U}^*$ should be understood as what the authors refer to as "sinks" in their framework. Similar analogies hold for the Consensus Process $\chi$. In particular, there are parallels with the notions introduced in [BG21] in the context of its continuous-time version. History graphs correspond to "oriented branching search paths" and their presence in $\chi^{\mathbf{1}}$ corresponds to being "consistent with respect to the initial configuration of the dominating process".

We now show that in a present history graph the cardinality of $\mathscr{U}^*$ cannot be too small, provided that either $j$ is sufficiently small or that the graph $G$ expands rapidly enough. Indeed, Lemma 4.5 and Lemma 4.6 below should be understood as giving a relationship between $j$ and the expansion constant of $G$ so that for every history graph $\mathscr{H}$, the cardinality of $\mathscr{U}^*$ is at least a positive proportion of the cardinality of $\mathscr{U}$. In the proof we will use the general identity

$$2|E(U)| = \sum_{x \in U} \deg(x) - |\partial_E U|, \quad (4.2)$$

where $E(U) := \{\{x,y\} \in E \mid x, y \in U\}$ denotes the set of edges of the subgraph of $G$ induced by the subset of vertices $U \subset V$.

**Lemma 4.5.** *Let $\Xi$ be either the $j$-threshold Consensus Process $\chi$ or the $j$-neighbour Bootstrap Percolation $\omega$ with $\varepsilon$-noise on $G$. Let $\mathscr{H} = (\mathscr{U}, \mathscr{F}, \mathscr{U}^*)$ be a history graph for $(o, T)$ for $\Xi$. Then*

$$(3\Phi_E - \Delta - 2j + 2)|\mathscr{U}| + \Delta - \Phi_E \leq 2(\Phi_E - j + 1)|\mathscr{U}^*|. \qquad (4.3)$$

*Moreover, if $\Xi = \chi$ and $G$ is bipartite, or if $\Xi = \omega$, we have the improvement*

$$2(\Phi_E - j + 1)|\mathscr{U}| + \delta - \Phi_E \leq (\Phi_E + \Delta - 2j + 2)|\mathscr{U}^*|. \qquad (4.4)$$

*Proof.* We first assume that $(o, T) \in \mathscr{U}^*$. In this case, Definition 4.1 forces $\mathscr{U} = \mathscr{U}^* = \{(o, T)\}$, since if $\mathscr{U}$ contained other vertices, item (1) would imply they all lie below $(o, T)$ and by item (3) this would make $(o, T)$ isolated in $\mathscr{H}$, a contradiction. Having settled this trivial case, we henceforth assume that $(o, T) \in \bar{\mathscr{U}} := \mathscr{U} \setminus \mathscr{U}^*$. We denote by $\mathscr{F}_t$ the edges in the history graph whose tail has time coordinate equal to $t$ and by

$$F_t := \{\{x, y\} \in E \mid \exists \, ((x, t), (y, t-1)) \in \mathscr{F}_t\}$$

the projection of $\mathscr{F}_t$ onto the base graph. We also recall that, according to the notation introduced at the beginning of Section 4, $W_t := \{x \in V \mid (x, t) \in \mathscr{W}\}$ denotes the space projection of the $t$-section of a generic $\mathscr{W} \subset \mathscr{V}$.

*Ineq. (4.5) for $\Xi = \chi$.* We observe that if $\{x, y\} \in F_t$, then either $x \in \bar{U}_t$ and $y \in U_{t-1}$, or vice versa, so that $F_t \subset E(\bar{U}_t \cup U_{t-1})$. It follows that

$$\left(E(\bar{U}_t) \sqcup \partial_E \bar{U}_t\right) \setminus E(\bar{U}_t \cup U_{t-1}) \subset \left(E(\bar{U}_t) \sqcup \partial_E \bar{U}_t\right) \setminus F_t = \bigcup_{x \in \bar{U}_t} \partial_E\{x\} \setminus F_t. \qquad (4.5)$$

On the one hand, the cardinality of the set on the right-hand side is upper bounded by $(j-1)|\bar{U}_t|$, because by the definition of history graph each vertex $(x, t) \in \mathscr{F}_t$ has exactly $\deg(x) - j + 1$ outgoing edges. On the other hand, we can lower bound twice the cardinality of the set on left-hand side by

$$2|E(\bar{U}_t) \sqcup \partial_E \bar{U}_t| - 2|E(\bar{U}_t \cup U_{t-1})| = |\partial_E \bar{U}_t| - \sum_{x \in U_{t-1} \setminus \bar{U}_t} \deg(x) + |\partial_E(\bar{U}_t \cup U_{t-1})|,$$

where we applied identity (4.1) twice and simplified the resulting sum over $\bar{U}_t$. Combining the two bounds and using the definition of edge expansion constant given by (2.3) for the set $\bar{U}_t \cup U_{t-1}$, we obtain:

$$|\partial_E \bar{U}_t| + \Phi_E |\bar{U}_t \cup U_{t-1}| - \sum_{x \in U_{t-1} \setminus \bar{U}_t} \deg(x) - 2(j-1)|\bar{U}_t| \leq 0, \qquad (4.6)$$

$$2(\Phi_E - j + 1)|\bar{U}_t| - (\Delta - \Phi_E)|U_{t-1} \setminus \bar{U}_t| \leq 0,$$

where for the second inequality we decomposed $\bar{U}_t \cup U_{t-1} = \bar{U}_t \sqcup (U_{t-1} \setminus \bar{U}_t)$ and used the estimates $\deg(x) \leq \Delta$ and $|\partial_E \bar{U}_t| \geq \Phi_E |\bar{U}_t|$. Now recall that $U_{t-1} = \bar{U}_{t-1} \sqcup U^*_{t-1}$. Summing over $t$ yields

$$2(\Phi_E - j + 1) \sum_{t=2}^{T} |\bar{U}_t| - (\Delta - \Phi_E) \sum_{t=2}^{T} |\bar{U}_{t-1}| \leq (\Delta - \Phi_E) \sum_{t=2}^{T} |U^*_{t-1}|,$$

$$(3\Phi_E - \Delta - 2j + 2)|\bar{\mathscr{U}}| + (\Delta - \Phi_E) \leq (\Delta - \Phi_E)|\mathscr{U}^*|,$$

where we used that $\bar{U}_1 = U^*_T = \emptyset$ and $\bar{U}_T = \{o\}$, because of item (1) and item (2) of Definition 4.1 and of the current assumption that $(o, T) \in \mathscr{U}$. This concludes the proof of (4.5) if $\Xi = \chi$. If $\Xi = \omega$, all the above arguments remain valid, but suboptimal, as we now show.

*Ineq. (4.5) for $\Xi = \omega$.* We proceed as in the general case above up to inequality (4.1). We then observe that, in the case of Bootstrap Percolation, the condition on straight edges in item (5) ensures that $\bar{U}_t \subset U_{t-1}$. Therefore, summing (4.1) over $t$ we obtain:

$$(\Phi_E - 2j + 2)\sum_{t=2}^{T}|\bar{U}_t| + \Phi_E \sum_{t=2}^{T}|U_{t-1}| - \sum_{t=2}^{T}\left(\sum_{x \in U_{t-1}} \deg(x) - \sum_{x \in \bar{U}_t} \deg(x)\right) \leq 0, \quad (4.7)$$

$$2(\Phi_E - j + 1)|\bar{\mathscr{U}}| + \Phi_E |\mathscr{U}^*| - \Phi_E - (\Delta|\mathscr{U}^*| - \deg(o)) \leq 0,$$

where in the second inequality we used the decomposition $U_{t-1} = \bar{U}_{t-1} \sqcup U^*_{t-1}$ and the identities $\bar{U}_1 = U^*_T = \emptyset$ and $\bar{U}_T = \{o\}$. The result (4.5) follows after rearranging and lower bounding $\deg(o) \geq \delta$.

*Ineq. (4.5) for $\Xi = \chi$ and $G$ bipartite.* We proceed as in the general case above up to inequality (4.1). We then observe that, in the case of Consensus Process on a bipartite graph, all vertices in $U_t$ are of the same color, opposite to that of $U_{t-1}$. In particular, this implies the identities $U_{t-1} \setminus \bar{U}_t = U_{t-1}$ and $|\partial_E \bar{U}_t| = \sum_{x \in \bar{U}_t} \deg(x)$, as $U_t$ is composed of isolated points. Therefore, summing (4.1) over $t$ we obtain exactly (4.1) and we conclude as in the case of Bootstrap Percolation. $\square$

**Lemma 4.6.** *Let $\Xi$ be either the $j$-threshold Consensus Process $\chi$ or the $j$-neighbour Bootstrap Percolation $\omega$ with $\varepsilon$-noise on $G$. Let $\mathscr{H} = (\mathscr{U}, \mathscr{F}, \mathscr{U}^*)$ be a history graph for $(o, T)$ for $\Xi$. If $\Xi = \chi$,*

$$(\Phi_V - j)|\mathscr{U}| + 1 \leq (\Phi_V - j)|\mathscr{U}^*|.$$

*If instead $\Xi = \omega$,*

$$(\Phi_V - j + 1)|\mathscr{U}| + 1 \leq (\Phi_V - j + 1)|\mathscr{U}^*|.$$

*Proof.* We sketch the proof only briefly, as the idea is very similar to the one of Lemma 4.5. In fact, considering the vertex boundary instead of the edge boundary simplifies the analysis. The analogue of (4.1) is

$$|\partial_V \bar{U}_t \setminus U_{t-1}| \leq \left|\bigcup_{x \in \bar{U}_t} \partial_V \{x\} \setminus U_{t-1}\right| \leq (j-1)|\bar{U}_t|.$$

If $\Xi = \chi$, we lower bound the left-hand side by $|\partial_V \bar{U}_t| - |U_{t-1}|$. If instead $\Xi = \omega$, using that in this case $\bar{U}_t \subset U_{t-1}$, we can improve the lower bound to $|\partial_V \bar{U}_t| + |\bar{U}_t| - |U_{t-1}|$. The conclusion follows from summing over all $t$, similarly to what is already done above. $\square$

*Remark* 4.7. Let $\omega$ be $j$-neighbour Bootstrap Percolation on $G$ (i.e. $\varepsilon = 0$). It is also natural to consider the maximal connected component $C_x$ of $\{y \in V \mid \lim_{t \to \infty} \omega_x^{\mu_q}(t) = 1\}$ containing a given $x \in V$ and ask if there exist $c > 0$ and $\tilde{q} > 0$ such that for all $q \in (0, \tilde{q}]$ and $\ell \geq 0$

$$\mathbb{P}[\text{diam}(C_x) \geq \ell] \leq q^{c(l+1)}.$$

If $G$ is the $d$-dimensional Euclidean lattice and $j > d$, this is [HT24b, Corollary 9.7] (in fact, their proof applies to all subcritical update rules). If instead $G$ is regular and $j > (d - \Phi_E)/2$, then the result follows from the proof of [BPP06, Theorem 1.4]. In this second case, an alternative, though more convoluted proof can be obtained by applying Proposition 5.2 to a decorated set system consisting either of "projections" of history graphs onto the base graph, or of shattered Toom countours, following the approach of [HS22].

## 4.2 Toom Argument

### 4.2.1 General setup

We begin by reformulating the definitions of Bootstrap Percolation and Consensus Process in the setting of [SST25b]. In order to do this, we need the following abstract definitions, see [SST25b, Section 2.2 and Definition 21].

**Definition 4.8** (Typed directed graph). A *typed directed graph* with *vertex set* $\mathsf{V}$, *vertex type set* $A$ and *edge type set* $B$ is a pair $(\mathcal{V}, \mathcal{E})$, where $\mathcal{V} \subset \mathsf{V} \times A$ and $\mathcal{E} \subset \mathsf{V} \times \mathsf{V} \times B$, such that for every $v \in \mathsf{V}$ there exists an $a \in A$ with $(v, a) \in \mathcal{V}$.

For each $a \in A$ and $b \in B$, we denote the set of vertices of type $a$ by $V_a := \{v \mid (v, a) \in \mathcal{V}\}$ and the set of directed edges of type $b$ by $\vec{E}_b := \{(v, w) \mid (v, w, b) \in \mathcal{E}\}$. Moreover, we call $(\mathsf{V}, \vec{E}) := (\bigcup_{a \in A} V_a, \bigcup_{b \in B} \vec{E}_b)$ the untyped directed graph associated to $(\mathcal{V}, \mathcal{E})$. The calligraphic font is dedicated to objects related to typed directed graphs.

**Definition 4.9** (Typed dependence graph and monotone cellular automata). A *typed dependence graph* with $\sigma \in \mathbb{N}_{\geq 1}$ types of edges is a typed directed graph $(\Lambda, \mathcal{H})$ with vertex type set $\{0, \bullet\}$ and edge type set $\{1, \ldots, \sigma\}$ such that vertices of type $0$ have no outgoing edges, vertices of type $\bullet$ have at least one outgoing edge of each type and the associated untyped directed graph is acyclic. The *monotone cellular automaton* associated with $(\Lambda, \mathcal{H})$ is the collection of maps $(\phi_i)_{i \in \Lambda}$ defined by

$$\phi_i(\xi) = \begin{cases} \sup_{s \in \{1,\ldots,\sigma\}} \inf_{\{j \mid (i,j) \in \vec{H}_s\}} \xi(j) & \text{if } i \in \Lambda_\bullet, \\ 0 & \text{if } i \in \Lambda_0. \end{cases} \quad (4.8)$$

We see that, given the initial configuration at time $0$, each given realization of the $j$-threshold Consensus Process with $\varepsilon$-noise is a monotone cellular automaton associated to the typed dependence graph defined by

$$\begin{aligned} \Lambda &= \Lambda_\bullet \sqcup \Lambda_0, \qquad \Lambda_\bullet = \{(x,t) \in \mathscr{V} \mid L_x(t) \geq \varepsilon\}, \quad \Lambda_0 = \{(x,t) \in \mathscr{V} \mid L_x(t) < \varepsilon\}, \\ \vec{H}_s &= \bigsqcup_{x \in V} \vec{H}_s^x, \qquad \vec{H}_s^x = \{((x,t),(y,t-1)) \mid y \in A_s^x, \, t \in \mathbb{N}_{\geq 1}\}, \end{aligned} \quad (4.9)$$

where, as $s$ varies, $A_s^x$ runs over the sets in the collection $\{A \mid A \subset N(x), |A| \geq j\}$ (recall that $N(x)$ denotes the set of neighbours of $x$ in the base graph $G$). The same is true for $j$-neighbour Bootstrap Percolation if we also allow $A_s^x$ to be $\{x\}$.

Having now recast our processes in the framework of [SST25b], we can apply their results to construct Toom cycles (when $\sigma = 2$) and countours (when $\sigma \geq 3$). A rough outline of the part of their program relevant in our context is as follows:

1. choose $\sigma \geq 2$ sets among the $(A_s^x)_s$, called *charges*, and consider the typed dependence (sub)graph they induce;

2. construct Toom cycles or contours with the special property of being present in the typed dependence graph of step 1;

3. show that the number of edges of a present Toom cycle or contour is bounded above by a (absolute) constant multiple of its sinks, which are the elements of a special subset of its vertices.

Regarding step 1, it is important to note that the framework of [SST25b] allows the charges to be chosen in a vertex-dependent way, without necessarily respecting any particular algebraic structure. To the best of our knowledge, this flexibility has not been exploited prior to the present work. Step 2, on the other hand, is the content of [SST25b, Section 4]. Since it is already designed to work on arbitrary graphs, it applies directly to our setting without modification. Step 3, instead, assumes that the underlying graph is a countable group. The reason is that it relies on so-called polar functions, which are homomorphisms that interact well with the charges chosen in step 1. Intuitively, charges represent directions in the base graph $G$ that are relevant for the dynamics of the process and the polar functions measure how much a given Toom cycle or contour expands along these directions. The set of sinks referred to in step 3, instead, plays a role analogous to that of the subset $\mathscr{U}^*$ of a history graph. As will become clear in Section 5, its cardinality reflects the probabilistic cost (or energy) associated with the presence of the corresponding Toom cycle or contour. For more precise formulations, we refer to [SST25b] for the case of $\mathbb{Z}^d$ and to [SST25a] for the case of arbitrary groups.

In the next two subsections we show how to carry out this program for the specific models considered in this work. As is common in such arguments, the main difficulty lies in finding a suitable model-specific polar function, since beyond the case of $\mathbb{Z}^d$, originally treated by Toom (see [SST25b, Lemma 10]), no general theory for this task is currently known.

### 4.2.2 Bootstrap Percolation

We consider $j$-neighbour Bootstrap Percolation with $\varepsilon$-noise on $G$ for $j \leq \bar{j}$, with $\bar{j}$ defined as in (2.3). Therefore, there exists $r \in V$ such that $j \leq |S_{k+1}(r) \cap N(x)|$ for every $k \in \mathbb{N}$ and $x \in S_k(r)$ (recall that $S_k(r)$ denotes the set of all vertices at distance exactly $k$ from $r$).

We define a typed dependence graph $(\Lambda, \mathcal{H})$ as in step 1 of the program described above by setting $\sigma = 2$ and defining, for each $x \in S_k(r)$, the charges

$$A_1^x := \{x\}, \qquad A_2^x \subset \{y \mid y \in S_{k+1}(r) \cap N(x)\} \text{ such that } |A_2^x| = j.$$

If the above condition does not uniquely determine $A_2^x$, as it happens in general, we choose it arbitrarily among the different possibilities, the condition $j \leq \bar{j}$ ensuring that there always exists at a least one valid option. In words, the untyped directed graph associated to $(\Lambda, \mathcal{H})$ is a subgraph of the space-time graph $\mathscr{G}$ of the previous subsection such that vertices corresponding to noise points have no outgoing edges, whereas vertices corresponding to update points have one outgoing straight edge, of type 1, and $j$ outgoing oblique edges, of type 2.

We now proceed with step 2. As $\sigma = 2$, we can specialize the construction of Toom contours to Toom cycles, which we now describe. For every fixed $n \in \mathbb{N}$, let $[n]$ denote the set $\{0, \ldots, n-1\}$ equipped with addition modulo $n$. We call the oriented graph $([n], \vec{E})$ an *oriented cycle* of length $n \geq 1$ if its set of oriented edges $\vec{E} \subset [n] \times [n]$ is such that for all $v \in [n]$ we have $|\vec{E} \cap \{(v, v+1), (v+1, v)\}| = 1$. The orientation chosen for the edges of an oriented cycle induces both a partition of its vertices, given by

$$V_\circ := \{v \in [n] \mid (v, v-1), (v, v+1) \in \vec{E}\}, \qquad V_* := \{v \in [n] \mid (v-1, v), (v+1, v) \in \vec{E}\},$$
$$V_1 := \{v \in [n] \mid (v-1, v), (v, v+1) \in \vec{E}\}, \qquad V_2 := \{v \in [n] \mid (v+1, v), (v, v-1) \in \vec{E}\},$$

and a partition of its edges, given by

$$\vec{E}_1 := \{(v, w) \in \vec{E} \mid w = v + 1\}, \qquad \vec{E}_2 := \{(v, w) \in \vec{E} \mid w = v - 1\}.$$

The elements of $V_*$ are called *sinks*. As a notational convenience, we might write $n$ instead of $0$ and, if $(0, n-1) \in \vec{E}_2$, we write $(n, n-1)$ instead of $(0, n-1)$.

The following is a rewriting of [SST25b, Definitions 24, 25].

**Definition 4.10.** A *Toom cycle in $\mathscr{V}$ rooted at $\psi_0$* is a triple $\mathcal{T} = ([n], \vec{E}, \psi)$, where $([n], \vec{E})$ is an oriented cycle of length $n$ such that $0 \in V_\circ$ and $\psi : \{0, \ldots, n\} \to \mathscr{V}$ is a map satisfying $\psi_0 = \psi_n$ and the following two conditions:

(1) for every $v \in V_*$ and $w \in [n]$, if $v \neq w$ then $\psi_v \neq \psi_w$;

(2) for every $v \in V_s \cup \{0\}$ and $w \in V_t \cup \{0\}$, with $s,t \in \{1, \circ, 2\}$, if $s \leq t$ and $\psi_v = \psi_w$ then $v \leq w$, where the order in $\{1, \circ, 2\}$ is defined by $1 < \circ < 2$.

Moreover, $\mathcal{T}$ is said to be *present* in the typed dependence graph $\mathcal{G}$ if $n \geq 2$ and

- $\psi_v$ is of type 0 for every $v \in V_*$;
- $(\psi_v, \psi_w)$ is of type $s$ for every $s \in \{1, 2\}$ and $(v, w) \in \vec{E}_s$ such that $v \in V_s \cup \{0\}$;
- $(\psi_v, \psi_w)$ is of type $3 - s$ for every $s \in \{1, 2\}$ and $(v, w) \in \vec{E}_s$ such that $v \in V_\circ \setminus \{0\}$.

If $n = 1$, only the first condition applies.

The following is a rewriting of [SST25b, Theorem 26] in our setting.

**Theorem 4.11.** *Let $\omega^1$ be $j$-neighbour Bootstrap Percolation with $\varepsilon$-noise on $G$, started from the constant one configuration. Assume that $\omega^1_o(T) = 0$ for some $(o, T) \in \mathscr{V}$. Then there exists a Toom cycle rooted at $(o, T)$ that is present in the typed dependence graph $\mathcal{G}$.*

We now proceed to step 3 of the program outlined in Section 4.2.1. For the specific model under consideration, it turns out that $|V^*|$ can in fact be computed explicitly in terms of the number of edges of the cycle. To this end, we define the function $L = (L_1, L_2) \colon \mathscr{V} \to \mathbb{Z}^2$ by setting $L_1(x, t) = -L_2(x, t)$ and $L_2(x, t) = \text{dist}(r, x)$, where dist denotes the graph distance on the base graph $G$. By construction, $L_1(x, t) + L_2(x, t) = 0$ for every $(x, t) \in \mathscr{V}$, and so $L$ is a *polar function* in the sense of [SST25b, Equation (5.18)]. In particular, it satisfies the zero-sum property

$$\sum_{s=1}^{2} \sum_{(v,w) \in \vec{E}_s} (L_s(\psi_w) - L_s(\psi_v)) = 0, \tag{4.10}$$

as proved in [SST25b, Lemma 42]. Alternatively, this identity can be verified directly in our setting by expanding the relevant definitions and rearranging terms.

The following version of [SST25b, Lemma 44], adapted to our setting, is what allows us to complete step 3. We include a full proof since, as explained above, their original argument relies on structural assumptions that do not apply here: in general, $G$ is not a group and $L$ is not a homomorphism. Furthermore, due to the simplicity of the update rule in our model, we are able to obtain an exact equality rather than an inequality.

**Lemma 4.12.** *Let $\mathcal{T} = ([n], \vec{E}, \psi)$ be a Toom cycle rooted at $(o, T)$ and present in $\mathcal{G}$. It holds that*

$$|V_\circ| = |V_*|, \qquad |\vec{E}_1| = |\vec{E}_2|, \qquad |\vec{E}| = 4(|V_*| - 1).$$

*Proof.* We refer to [SST25b, Section 5.2] for the first two equalities, as our argument follows theirs without modification. It is also possible to verify them directly with little difficulty. Regarding the last one, by first using the definition of $\vec{E}_1$ and $\vec{E}_2$ and then the one for typed vertices, we can decompose the zero sum property (4.2.2) as follows:

$$0 = \sum_{\substack{v \in V_1 \cup \{0\} \\ (v, v+1) \in \vec{E}}} L_1(\psi_{v+1}) - L_1(\psi_v) + \sum_{\substack{v \in V_\circ \setminus \{0\} \\ (v, v+1) \in \vec{E}}} L_1(\psi_{v+1}) - L_1(\psi_v)$$
$$+ \sum_{\substack{v \in V_2 \cup \{0\} \\ (v, v-1) \in \vec{E}}} L_2(\psi_{v-1}) - L_2(\psi_v) + \sum_{\substack{v \in V_\circ \setminus \{0\} \\ (v, v-1) \in \vec{E}}} L_2(\psi_{v-1}) - L_2(\psi_v)$$

where the conditions $(v, v + 1) \in \vec{E}$ and $(v, v - 1) \in \vec{E}$ are actually redudant, as they are determined by the vertex type. We now observe that the first sum in the above display vanishes, as for each vertex $v$ in that sum, the edge $(\psi_{v+1}, \psi_v)$ is straight. Similarly, the fourth sum also vanishes. In contrast, each term in the second sum evaluates to $-1$, while each term in the third sum equals 1. Therefore, we have

$$0 = -|V_\circ \setminus \{0\}| + |V_2 \cup \{0\}| = -(|V_*| - 1) + |\vec{E}_2| - (|V_*| - 1),$$

from which the conclusion follows. $\square$

### 4.2.3 Toom contours for Consensus Process on regular trees

We consider the $j$-threshold Consensus Process with $\varepsilon$-noise $\chi$ on the $d$-regular tree $\mathbb{T}_d$ with $j < d$ as the monotone cellular automaton associated to the typed dependence graph defined by (4.2.1). We now follow the plan outlined in Section 4.2.1 above.

We begin step 1 by orienting the edges of $\mathbb{T}_d$ so that each vertex has exactly $d-1$ incoming edges and one outgoing edge. This can be achieved by first selecting a simple infinite oriented path, and then orienting all remaining edges toward this path. For every vertex $x$, we denote its $d$ neighbours by $y_i^x$, with $i \in \{1, \ldots, d\}$, in such a way that the oriented edges are $(y_i^x, x)$ for $i \in \{1, \ldots, d-1\}$ and $(x, y_d^x)$. We then choose the $\sigma = d$ charges defined by

$$A_s^x := \{y_{s'}^x \mid s' \in \{1, \ldots, d\} \setminus \{s\}\}, \qquad s \in \{1, \ldots, d\}, \quad x \in V. \tag{4.11}$$

In particular, we observe that the dependence on $j$ has disappeared. Effectively, we have reduced ourselves to the case $j = d - 1$, for which establishing Proposition 2.9 is most challenging, since for any $j' \leq j$ the $j'$-threshold Consensus Process stochastically dominates the one with threshold $j$.

Step 2 consists in the application of [SST25b, Theorem 23], which in our setting states that, if $\chi_o^{\mathbf{1}}(T) = 0$ for some $(o, T) \in \mathscr{V}$, then there exists a Toom contour rooted at $(o, T)$ that is present in the typed dependence graph induced by (4.2.3). Toom contours are a generalization of the Toom cycles defined in Section 4.2.2. In particular, they also come with a set of sinks $V_*$ and allow for an arbitrary (finite) number $S$ of vertex and edge types, as opposed to just two. The downside is that precise definitions and statements are technically involved. For this reason we do not repeat them here, but refer to [SST25b, Definitions 16, 17, 18 20, 21, 22 and Theorem 23].

We finally come to step 3. Let $\bar{x} \in V$ be an arbitrarily chosen vertex. We define $L := (L_1, \ldots, L_d) \colon V \to \mathbb{Z}^d$ by setting $L(\bar{x}) = 0$ and requiring that, for every vertex $x \in V$, the following relations hold:

$$L(y_i^x) - L(x) = -e_i + e_d, \qquad i \in \{1, \ldots, d-1\},$$

where $e_i$ denotes the standard basis vector in $\mathbb{Z}^d$ with one in the $i$-th coordinate and zero elsewhere. This definition is well-posed, as the above equations are a system of linear recurrence relations in $L$ indexed by the vertices of $\mathbb{T}_d$ and with initial value 0 at $\bar{x}$. Morevoer, it follows by induction that

$$L_d(y_d^x) - L_d(x) = -1, \qquad \sum_{i=1}^d L_i(x) = 0. \tag{4.12}$$

We trivially lift $L$ to the space-time vertices $\mathscr{V}$ by setting

$$L_s(x, t) := L_s(x), \qquad (x, t) \in \mathscr{V}, \quad s \in \{1, \ldots, d\}.$$

The last equation in (4.2.3) implies that $L$ is a *polar function* in the sense of [SST25b, Equation (5.18)]. In particular, it satisfies the zero-sum property of [SST25b, Lemma 42], which is the analogue of (4.2.2) for an arbitrary finite number of types. Finally, analogously to [SST25b, Definition 6], we define the *edge speed* of $\phi_{x,t}$ in the direction $L_s$ to be

$$\varepsilon_s(\phi_{x,t}) := \sup_{s' \in \{1, \ldots, d\}} \inf_{y \in A_{s'}^x} \{L_s(y, t-1) - L_s(x, t)\}, \qquad s \in \{1, \ldots, d\},$$

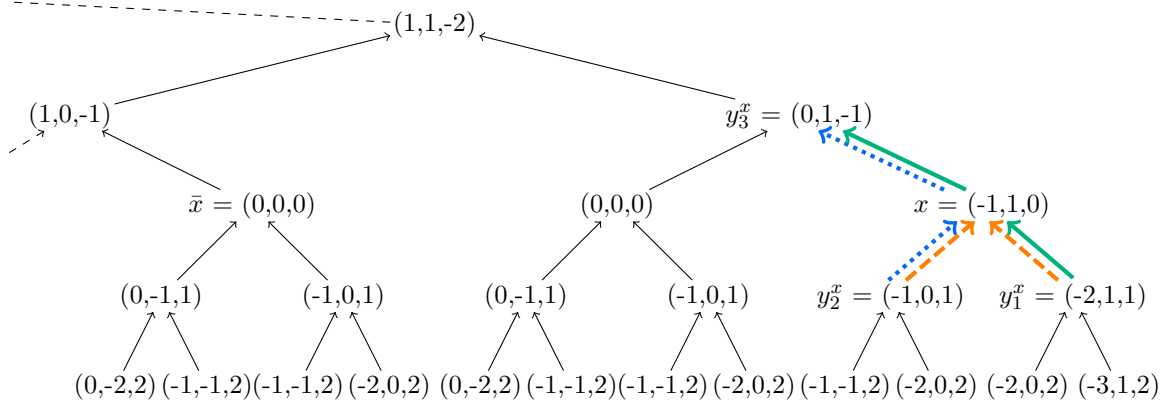

Figure 1: A portion of the 3-regular tree $\mathbb{T}_3$, oriented so that each vertex, labeled by the value assigned to it by the polar function $L$, has two incoming edges and one outgoing edge. A specific vertex $x$ is highlighted, along with its neighbours $y_1^x$, $y_2^x$, and $y_3^x$. Dotted blue, thick green, and dashed orange edges connect $x$ to the vertices in the sets $A_1^x$, $A_2^x$, and $A_3^x$, respectively. The embedding in the plane is chosen so that each level set of $L_3$ lies on a distinct horizontal line. On any such line, pairs of vertices sharing the same parent are drawn in decreasing order of $L_1$ (equivalently, increasing order of $L_2$).

where we recall that $\phi_{x,t}\colon \{0,1\}^{\mathscr{V}} \to \{0,1\}$, defined as in (4.9), is the local function which determines the value of $\chi_x(t)$ based on the configuration of the process at time $t-1$. By direct computation, we obtain that

$$\varepsilon_s(\phi_{x,t}) = 0, \qquad \varepsilon_d(\phi_{x,t}) = 1, \qquad (x,t) \in \mathscr{V}, \qquad s \in \{1,\ldots,d-1\},$$

$$R_s := -\inf\left\{ L_s(y, t-1) - L_s(x, t) \,\bigg|\, y \in \bigcup_{s'=1}^d A_{s'}^x \right\} = 1, \qquad s \in \{1,\ldots,d\}.$$

In particular, we observe that neither $\varepsilon_s$ nor $R_s$ depend on $(x,t)$.

With the current notation, the bound provided by [SST25b, Lemma 41] or, equivalently, by [SST25a, Lemma 5.2], takes the form

$$n_E(\mathcal{T}) \leq 3\left(1 + \frac{\sum_{s=1}^d R_s}{\sum_{s=1}^d \varepsilon_s}\right)(|V_*| - 1) = 3(d+1)(|V_*| - 1), \qquad (4.13)$$

where $n_E(\mathcal{T})$ denotes the number of the edges of the present Toom contour $\mathcal{T}$ rooted at $(o, T)$ and $V_*$ is the set of its sinks. Although the proofs of these lemmas assume that the polar function $L$ is either linear or a homomorphism, assumptions which do not hold in our setting[1], what is really required is that $L$ satisfies $\sum_{s=1}^d \varepsilon_s > 0$ and $\max\{R_s \mid s \in \{1,\ldots,d\}\} < \infty$, uniformly over all space-time vertices. As these conditions are satisfied by our construction, step 3 is also achieved.

## 5 Decorated Set Systems

We begin this section by recalling the abstract framework of decorated set systems introduced in [HT24b], along with the associated result establishing exponential decay properties. We then proceed to construct different instances of decorated set systems, which will be used in the proof of Proposition 2.9.

---

[1] It is possible to describe the graph $\mathbb{T}_3 \square \mathbb{Z}$ as a Cayley graph, for instance using the presentation $\langle a, b, d, d^{-1} \mid b^2, [a,d], [b,d]\rangle \cong \langle a, b \mid b^2 \rangle \times \langle d, d^{-1}\rangle$, where $[a,d] := ada^{-1}d^{-1}$ denotes the commutator of the two elements $a$ and $d$ and $\times$ denotes the direct product of groups. However, we were not able to exploit this to obtain homomorphisms that also serve as polar functions with nontrivial edge speeds, as required in [SST25a].

Decorated set systems can, in principle, be defined on any graph, though the setting of interest is the space-time graph associated with the process under consideration. In the original work [HT24b, Section 9], this was the lattice $\mathbb{Z}^D$, with $D = d+1$ and $d$ denotes the spacial dimension. The analogue in our case would be the Cartesian product graph $G \square \mathbb{N}$. However, we will work instead with the strong product graph $G \boxtimes \mathbb{N}$, as this choice is both more natural in our context and yields stronger results, as explained at the end of Section 2. Throughout, the notions of diameter and path are thus understood with respect to $G \boxtimes \mathbb{N}$. For convenience, we set $\mathscr{V} \coloneqq V \times \mathbb{N}$.

**Definition 5.1** (Decorated set system). For any $\mathscr{Z} \subset \mathscr{V}$, let $\Gamma_{\mathscr{Z}}$ be a (possibly empty) arbitrary set. We call $(\mathscr{Z}, \gamma)$ a *decorated set* if $\mathscr{Z}$ is nonempty and bounded and $\gamma \in \Gamma_{\mathscr{Z}}$. Two decorated sets $(\mathscr{Z}_1, \gamma_1)$ and $(\mathscr{Z}_2, \gamma_2)$ are *disjoint* if $\mathscr{Z}_1 \cap \mathscr{Z}_2 = \emptyset$. A *decorated set system* is a tuple $(\mathsf{P}, \mathsf{E})$, where $\mathsf{P}$ is a probability measure and $\mathsf{E}$ is a function which associates to each decorated set $(\mathscr{Z}, \gamma)$ an event $\mathsf{E}(\mathscr{Z}, \gamma)$ in such a way that for any finite collection of disjoint decorated sets $(\mathscr{Z}_i, \gamma_i)_{i \in I}$ we have

$$\mathsf{P}\left[\bigcup_{i \in I} \mathsf{E}(\mathscr{Z}_i, \gamma_i)\right] \leq \prod_{i \in I} \mathsf{P}\left[\mathsf{E}(\mathscr{Z}_i, \gamma_i)\right]. \tag{5.1}$$

**Proposition 5.2.** *Consider a decorated set system. Assume that for some $C > 1$ and $\epsilon > 0$ small enough depending on $C$ the following hold:*

1. *for all decorated sets $(\mathscr{Z}, \gamma)$ we have $\mathrm{diam}(\mathscr{Z}) \leq C|\mathscr{Z}|$;*
2. *for every $(x,t) \in \mathscr{V}$ the number of decorated sets $(\mathscr{Z}, \gamma)$ such that $|\mathscr{Z}| = n$ and $(x,t) \in \mathscr{Z}$ is at most $C^n$;*
3. *for all decorated sets $(\mathscr{Z}, \gamma)$ we have $\mathsf{P}[\mathsf{E}(\mathscr{Z}, \gamma)] \leq \epsilon^{|\mathscr{Z}|/C}$.*

*Fix $n \in \mathbb{N}_{\geq 1}$ and $(x,t) \in \mathscr{V}$. Then the probability that $(x,t)$ is the endpoint of a path $\mathscr{P}$ in $G \boxtimes \mathbb{N}$ of length at least $n$ and such that for every $p \in \mathscr{P}$ there exists a decorated set $(\mathscr{Z}_p, \gamma_p)$ such that $p \in \mathscr{Z}_p$ and $\mathsf{E}(\mathscr{Z}_p, \gamma_p)$ occurs is at most $\epsilon^{n/(7C+7)^2}$.*

Proposition 5.2 is proved in [HT24b] for the case where the underlying graph is $\mathbb{Z}^D$. However, their proof does not rely on any specific property of the Euclidean lattice and is already phrased entirely in terms of graph distance, so that it applies directly to our setting. The only exeption is their inequality (9.7), which in general does not hold in our context, since the cardinality of a ball in $G \boxtimes \mathbb{N}$ may grow exponentially with its radius. Nevertheless, due to our assumption that the degrees of $G$ are uniformly bounded, the upper bound provided by that inequality can be replaced by a quantity growing at most exponentially in the total number of vertices in the chain (a chain is a sequence of decorated sets, satisfying some useful properties, which is used in their proof). This is sufficient to ensure that the rest of their argument carries over without further modifications, provided the $\epsilon$ of Proposition 5.2 is chosen sufficiently small, depending on the maximum degree $\Delta$ of $G$.

**Lemma 5.3.** *Let $\Xi$ be either the $j$-threshold Consensus Process $\chi$ or the $j$-neighbour Bootstrap Percolation $\omega$ with $\varepsilon$-noise on $G$. Define the function $\mathsf{E}$ on the collection of decorated sets*

$$\{\mathscr{H} = (\mathscr{U}, \mathscr{F}) \mid \mathscr{H} \text{ is a history graph for } (o, T) \text{ for } \Xi \text{ for some } (o, T) \in \mathscr{V}\},$$

*where $\mathscr{U}$ is the set and $\mathscr{F}$ is the decoration, by setting*

$$\mathsf{E}(\mathscr{U}, \mathscr{F}) = \{\zeta \in \{0,1\}^{\mathscr{V}} \mid L_x(t) < \varepsilon \text{ for every } x \in \mathscr{U}^*\}.$$

*Assume that $j$ is $\Xi$-good as in any of the items (a), (i), (ii) and (iii) of Definition 2.5. Then $(\mathbb{P}, \mathsf{E})$ is a decorated set system in $G \boxtimes \mathbb{N}$ satisfying the hypotheses of Proposition 5.2.*

*Proof.* First of all, observe that Remark 4.2 ensures that $\mathscr{X}$ and $\mathsf{E}$ as above are well-defined. It is clear that, when measuring the probability of the events $\mathsf{E}(\mathscr{U}, \mathscr{F})$, $\mathbb{P}$ acts like a Bernoulli product measure. Therefore, if $\mathscr{U}_1 \cap \mathscr{U}_2 = \emptyset$, then the events $\mathsf{E}(\mathscr{U}_1, \mathscr{F})$ and $\mathsf{E}(\mathscr{U}_2, \mathscr{F})$ are independent under $\mathbb{P}$, so that in particular (5.1) holds. This proves that $(\mathbb{P}, \mathsf{E})$ is a decorated set system.

We now verify that $(\mathbb{P}, \mathsf{E})$ satisfies the conditions of Proposition 5.2. Condition 1 clearly holds for any $C \geq 1$, since every history graph $\mathscr{H}$ is connected in $\mathscr{G}$ and thus, *a fortiori*, also in $G \boxtimes \mathbb{N}$.

To verify condition 2, we first recall that every history graph $\mathscr{H} = (\mathscr{U}, \mathscr{F})$ is uniquely determined by its edge set $\mathscr{F}$ (see Remark 4.2). Furthermore, we recall from Section 4 the decompositions $\mathscr{U} =: \bar{\mathscr{U}} \sqcup \mathscr{U}^*$ and $U_t = \bar{U}_t \sqcup U_t^*$ for every $t \in [0,T]$, where $U_t$ is the space projection of the $t$-section of $\mathscr{U}$. We then observe that any history graph for a fixed $(o,T) \in \mathscr{V}$ can be constructed as in the algorithm described in the proof of Lemma 4.3. In particular, at each time $t \in [1,T]$, the generic vertex $x \in U_t$ can either belong to $\bar{U}_t$ or to $U_t^*$ and, if it belongs to $\bar{U}_t$, then step (i) of the algorithm gives at most $\binom{\deg(x)}{\deg(x)-j+1} \leq 2^\Delta$ possible choices to select its outgoing edges. Now suppose that $|\mathscr{U}| = n$. Then $U_t \neq \emptyset$ implies $t \in [\max\{T-n+1,1\}, T]$, because $\mathscr{H}$ is connected in $G \boxtimes \mathbb{N}$. In the following estimate, we will use the symbols $m_1, \ldots, m_n \in \mathbb{N}$ to represent the cardinalities of $U_T, \ldots, U_{\max\{T-n+1,0\}}$, purposefully ignoring that, in fact, a history graph might have less than $n$ nonempty $t$-sections, $m_1$ is always equal to 1 and $m_i$ can only be zero if it is in the tail of the sequence $(m_1, \ldots, m_n)$. The only constraint we use is $\sum_{i=1}^n m_i = n$. Decomposing over the possible cardinalities of the $U_t$'s, we obtain the following upper bound for the number of history graphs for a fixed $(o,T)$ with vertex set of cardinality $n$:

$$\sum_{m_1+\cdots+m_n=n} \prod_{i=1}^n \left(1+2^\Delta\right)^{m_i} = \binom{2n-1}{n}\left(1+2^\Delta\right)^n \leq 2^{2n}\left(1+2^\Delta\right)^n.$$

Next, consider a history graph $\mathscr{H} = (\mathscr{U}, \mathscr{F})$ for some $(o,T)$ with $|\mathscr{U}| = n$ and containing a fixed vertex $(x,t) \in \mathscr{U}$. Since $\mathscr{U}$ is connected in $G \boxtimes \mathbb{N}$, the vertex $(o,T)$ must belong to $\mathscr{B}_n((x,t))$, the ball in $G \boxtimes \mathbb{N}$ of radius $n$ centered at $(x,t)$. As the maximum degree $\Delta$ of $G$ is finite, we have $|B_n((x,t))| \leq (1+\Delta)^n$ and $|\mathscr{B}_n((x,t))| \leq (1+2n)(1+\Delta)^n$. Putting everything together, condition 2 is satisfied for any $C$ such that $C^n \geq (1+2n)(1+\Delta)^n 4^n (1+2^\Delta)^n$.

Finally, condition 3 follows directly from either Lemma 4.5 or Lemma 4.6, depending on which of the items of Definition 2.5 is satisfied by $j$. $\square$

**Lemma 5.4.** *Let $\omega$ be $j$-neighbour Bootstrap Percolation with $\varepsilon$-noise on $G$. Define the function $\mathsf{E}$ on the collection of decorated sets*

$$\{\mathcal{T} = ([n], \vec{E}, \psi) \mid \mathcal{T} \text{ is a present Toom cycle rooted at } (o,T) \text{ for some } (o,T) \in \mathscr{V}\},$$

*where the set is $\psi([n])$ and the decoration consists of both the edge set $\vec{E}$ and the map $\psi$, by setting*

$$\mathsf{E}(\mathcal{T}) = \{\zeta \in \{0,1\}^{\mathscr{V}} \mid L_x(t) < \varepsilon \text{ for every } x \in V_*\}.$$

*Assume that $j$ is $\omega$-good as in item (b) of Definition 2.5. Then $(\mathbb{P}, \mathsf{E})$ is a decorated set system in $G \boxtimes \mathbb{N}$ satisfying the hypotheses of Proposition 5.2.*

Given the result of Lemma 4.12 and in light of [SST25b, Lemma 43], which establishes that the number of present Toom cycles grows at most exponentially in the number of their edges, the proof proceeds analogously to that of Lemma 5.3 and is therefore omitted.

# 6 Proofs of the main results

*Proof of Proposition 2.9.* The proof follows directly by applying Proposition 5.2 to different decorated set systems, depending on which items of Definition 2.5 hold.

If we are considering the $j$-neighbour Bootstrap Percolation and item (b) holds, we apply Theorem 4.11 and Lemma 5.4, which rely on Toom cycles. If instead we are considering the $j$-threshold Consensus Process on a $d$-regular tree with $j < d$ (i.e. item (iv) holds), we use the Toom contours described in Section 4.2.3. The analogue of Theorem 4.11 is then [SST25b, Theorem 23], while the analogue of Lemma 5.4 is obtained by verifying the conditions of Proposition 5.2: condition 1 follows immediately by connectedness, condition 2 is established in [SST25b, Lemma 40] and (4.2.3) ensures condition 3. In all remaining cases, we apply Lemma 4.3 and Lemma 5.3. □

*Proof of Theorem 2.6.* Let $(\eta^{\mu_p}(t))_{t\geq 0}$ be the FA-$j$f model on $G$ with initial distribution $\mu_p = \text{Ber}(p)^{\otimes V}$, equilibrium parameter $q$ and $\omega$-good $j$. We consider the initial configurations of the two processes $\eta^{\mu_p}$ and $\eta^{\mu_q}$ to be monotonically coupled. That is, almost surely, if $p \leq q$, then $\eta_x^{\mu_p}(0) \leq \eta_x^{\mu_q}(0)$ for all $x \in V$, whereas the reverse inequality holds if $p > q$. Since both processes evolve according to the same realization of the Poisson clocks and the same family of random variables $(U_x(t))_{x,t}$, this yields a coupling of the entire processes over time. Given that $\eta^{\mu_q}$ is stationary, the left-hand side of (2.6) is bounded above by

$$|\mathbb{E}[f(\eta^{\mu_p}(t)) - \mu_q(f)]| \leq \|f\|_\infty \, \mathbb{P}[\eta_S^{\mu_p}(t) \neq \eta_S^{\mu_q}(t)] \leq \|f\|_\infty |S| \, \mathbb{P}[\eta_o^{\mu_p}(t) \neq \eta_o^{\mu_q}(t)],$$

where $\eta_S(t)$ denotes the restriction of the configuration $\eta(t)$ to the vertices in $S := \text{supp}(f)$ (which is a finite set, $f$ being a local function) and $o \in S$ is a vertex maximizing the probability of the event $\{\eta_x^{\mu_p}(t) \neq \eta_x^{\mu_q}(t)\}$ over all vertices $x \in S$.

We now claim that, for the $\widetilde{\varepsilon}$ of Proposition 2.9, there exists $T = T(\widetilde{\varepsilon}) > 0$ sufficiently large and $\varepsilon = \varepsilon(\widetilde{\varepsilon}, T)$ sufficiently small such that for all $p, q \in [1 - \varepsilon, 1]$,

$$\mathbb{P}[\eta_o^{\mu_p}(t) \neq \eta_o^{\mu_q}(t)] \leq \mathbb{P}[\forall (y, s) \in \mathcal{K}_{x,t} \; \omega_y^{\mathbf{1}}(s) = 0] \leq \varepsilon_0^{c \lfloor t/T \rfloor}, \qquad (6.1)$$

where $\omega^{\mathbf{1}}$ is $j$-neighbour Bootstrap Percolation with $\widetilde{\varepsilon}$ noise started from the constant one configuration and $\mathcal{K}_{x,t}$ is a connected component of space-time points containing both $(o, \lfloor t/T \rfloor)$ and a vertex in $V \times \{1\}$. The second inequality follows immediately from Proposition 2.9. The first one, instead, is derived in [HT24b, Equations (10.1)-(10.5)] for the case where $G$ is the $d$-dimensional Euclidean lattice. It relies on coupling $\eta$ with $\omega$ via a version of the Contact Process, together with a renormalization both in space and time. It turns out that the only part of their arguments which uses properties of the Euclidean lattice is the renormalization in space, which in our setting is not necessary. This is due to the fact that the functions $c_x$ in (2.1) only consider the nearest neighbours of $x$, unlike the more general dynamics considered in [HT24b]. Therefore, we do not repeat the proof and instead ouline the main ideas while providing precise references. We point out that a simplified rewriting for the specific case of FA-2f on $\mathbb{Z}^2$ can also be found in [HT24a, Section 7.3].

First, we define the Contact Process $\zeta$ with parameter $q_0 \in [0, 1]$ by modifying Definition 2.1 so that $\zeta_x(t) = 0$ whenever $U_x(t) > q_0$, regardless of the value of $c_x(\zeta(t^-))$. By [HT24b, Corollary 8.2], if the initial configuration $\zeta(0)$ is sampled from a Bernoulli product measure with a sufficiently high density, it is possible to renormalize time into intervals of length $T$ so that, upon choosing the parameter $q_0$ large enough, one obtains the stochastic domination $\omega^{\mathbf{1}} \lesssim \widehat{\omega}$, where $\widehat{\omega}$ is defined by

$$\widehat{\omega}_x(\tau) := \mathbb{1}_{\{\forall (y,t) \in \{x\} \times [T(\tau-1), T\tau), \, \zeta_y(t) = 1\}} \qquad \forall (x, \tau) \in V \times \mathbb{N}_{\geq 1}. \qquad (6.2)$$

We remark that this makes use of the Liggett-Schonmann-Stacey theorem [LSS97; Swa22], which is necessary because the "good boxes" used in [HT24b, Lemma 8.1], upon which [HT24b,

Corollary 8.2] relies, exhibit nontrivial (but bounded) time dependencies. In our context, since we do not renormalize in space, we effectively set their constant $R = 1$.

Second, as shown in [HT24b, Lemma 7.1], the processes $\eta^{\mu_p}$ and $\eta^{\mu_q}$ can be coupled provided that their equilibrium parameter $q$ is larger than $q_0$ and that their initial distributions $\mu_p$ and $\mu_q$ stochastically dominate the initial distribution of $\zeta$. By coupled we mean that, at each time $t$, the two configurations $\eta^{\mu_p}(t)$ and $\eta^{\mu_q}(t)$ can differ only on the random subset of "orange healthy sites" $O_t \subset \{(y,t) \in V \times [0,\infty) \,|\, \zeta_y(t) = 0\}$, as defined in [HT24b, Equation (7.2)]. The connected set $\mathscr{K}_{(x,t)}$ in (6) is obtained through (6) as a renormalized, discrete-time version of the space-time cluster $\bigcup_{s \geq 0} O_s$. $\square$

*Proof of Theorem 2.7.* Let $\chi$ be the $j$-threshold Consensus Process with $\varepsilon$-noise on $G$ and with $\chi$-good $j$ satisfying (i), (ii) or (iii) of Definition 2.5. By Lemma 4.3, if $\chi_o^{\mathbf{1}}(T) = 0$, then there exists a history graph $\mathscr{H} = (\mathscr{U}, \mathscr{F})$ for $(o, T)$ which is present in $\chi^{\mathbf{1}}$. We denote this event by $\mathsf{E}(\mathscr{U}, \mathscr{F})$ and note that

$$\mathbb{P}[\mathsf{E}(\mathscr{U}, \mathscr{F})] = \varepsilon^{|\mathscr{U}^*|} \leq \varepsilon^{|\mathscr{U}|/K}, \tag{6.3}$$

where $K > 0$ is a constant, independent of $\mathscr{H}$, provided by either Lemma 4.5 or Lemma 4.6, depending on which item of Definition 2.5 we are considering. First combining these observations, then decomposing over the cardinality of $\mathscr{U}$ and finally applying a union bound, we obtain:

$$\mathbb{P}[\chi_o^{\mathbf{1}}(T) = 0] \leq \mathbb{P}[\text{a history graph } \mathscr{H} = (\mathscr{U}, \mathscr{F}) \text{ for } (o,T) \text{ is present in } \chi^{\mathbf{1}}]$$
$$\leq \sum_{n=1}^{\infty} \sum_{\mathscr{H}} \mathbb{P}[\mathsf{E}(\mathscr{U}, \mathscr{F}), |\mathscr{U}| = n] \leq \sum_{n=1}^{\infty} C^n \varepsilon^{n/K}, \tag{6.4}$$

where the last inequality follows from the same argument used in the proof of Lemma 5.3, which shows that the number of history graphs $\mathscr{H}$ for $(o, T)$ with $|\mathscr{U}| = n$ is at most $C^n$ for some (explicit) constant $C > 1$. We observe that the bound provided by (6) is uniform in $(o, T)$ and that the series on the right-hand side is less than $1/2$ for sufficiently small $\varepsilon$. Since $\chi^{\mathbf{0}}(t) = \mathbf{0}$ for all $t \geq 0$, it follows that the process $\chi$ admits more than one invariant measure when $\varepsilon$ is sufficiently small.

We now consider the $\varepsilon$-noise Majority Vote Process $\sigma$ on $G$. From the definitions of Section 2, it follows immediately that $\sigma$ is monotonically coupled with the $j := \lfloor \Delta/2 + 1 \rfloor$-threshold Consensus Process $\chi$, so that for any initial configuration $\xi \in \Omega$, time $t \geq 1$ and vertex $x$ we have $\mathbb{P}[\chi_x^\xi(t) \leq \sigma_x^\xi(t)] = 1$. In particular, $\mathbb{P}[\sigma_o^{\mathbf{1}}(T) = 0] \leq \mathbb{P}[\chi_o^{\mathbf{1}}(T) = 0]$. Assuming that $\lfloor \Delta/2 + 1 \rfloor$ is $\chi$-good, it then follows from (6) that $\mathbb{P}[\sigma_o^{\mathbf{1}}(T) = 0] \leq c$ for some $c \in (0, 1/2)$. Therefore, by symmetry between the states 0 and 1, we conclude that the process $\sigma$ admits more than one invariant measure.

Finally, let $d \geq 2$ and consider the $j$-threshold Consensus Process on the $d$-regular tree $\mathbb{T}_d$ with $j < d$, so that item (iv) of Definition 2.5 holds. We claim that the Peierls argument given in (6) carries over directly, with history graphs replaced by the Toom contours discussed in Section 4.2.3. Indeed, the presence of a Toom contour on the event $\chi_o^{\mathbf{1}}(T) = 0$ is guaranteed by [SST25b, Theorem 23], while the analogue of (6) is ensured by (4.2.3). Moreover, [SST25b, Lemma 40] establishes that the number of present Toom contours grows at most exponentially in the number of their edges. The result then follows from the same reasoning as above. $\square$

## Acknowledgements

The author would like to thank Ivailo Hartarsky for his invaluable guidance, constant support and numerous insightful discussions, as well as Réka Szabó and Fabio Toninelli for their helpful comments. The author also gratefully acknowledges financial support of the Austrian Science Fund (FWF), Project Number P 35428-N.